\documentclass[12pt]{article}
\usepackage{graphicx}
\usepackage{amsmath, amsthm}
\usepackage{amssymb}
\usepackage{amsfonts}

\newtheoremstyle{theorem}
  {10pt}		  
  {10pt}  
  {\sl}  
  {\parindent}     
  {\bf}  
  {. }    
  { }    
  {}     
\theoremstyle{theorem}

\newtheoremstyle{defi}
  {10pt}		  
  {10pt}  
  {\rm}  
  {\parindent}     
  {\bf}  
  {. }    
  { }    
  {}     
\theoremstyle{defi}

\begin{document}

\title{Central Configurations in the Trapezoidal Four-Body Problems}
\author{Muhammad Shoaib \\
Faculty of Sciences, Department of Mathematics,\\ University of Hail,
Kingdom of Saudi Arabia,\\email: safridi@gmail.com}
\maketitle

\begin{abstract}
In this paper we discuss the central configurations of the Trapezoidal
four-body Problem. We consider four point masses on the vertices of an
isosceles trapezoid with two equal masses $m_1=m_4$ at positions $(\mp 0.5,
r_B)$ and $m_2=m_3$ at positions $(\mp \alpha/2, r_A)$. We derive, both
analytically and numerically, regions of central configurations in the phase
space where it is possible to choose positive masses. It is also shown that
in the compliment of these regions no central configurations are possible.


\end{abstract}

\section{Introduction}

The classical equation of motion for the n-body problem has the form (\cite{shoaib1}- \cite{shoaib5})
\begin{equation}
m_{i}\frac{d^{2}\vec{r}_{i}}{dt^{2}}=\frac{\partial U}{\partial \vec{r}_{i}}%
=\sum_{j\neq i}\frac{m_{i}m_{j}\left( \vec{r}_{i}-\vec{r}_{j}\right) }{|\vec{%
r}_{i}-\vec{r}_{j}|^{3}}\qquad i=1,2,...,n,  \label{eqtnOfMotion-General}
\end{equation}%
where the units are chosen so that the gravitational constant is equal to
one, $\mathbf{r}_{i}$ is the location vector of the $i$th body,
\begin{equation}
U=\sum_{1\preceq i<j\preceq n}\frac{m_{i}m_{j}}{|\vec{r}_{i}-\vec{r}_{j}|}
\end{equation}%
is the self-potential, and $m_{i}$ is the mass of the $i$th body. To understand the dynamics presented by a total collision of the masses or
the equilibrium state of a rotating system, we are led to the concept of a central configuration (\cite{AlanAlbouy2002}- \cite{MK},\cite{ShoaibAIP} and \cite{ShoaibHindawi}). A \emph{central configuration} is a particular configuration of the $n$-bodies where the acceleration vector of each body is proportional to its position vector, and the constant of proportionality is the same for the $n$%
-bodies, therefore
\begin{equation}
\sum_{j=1,j\neq i}^{n}\frac{m_{j}(\vec{r}_{j}-\vec{r}_{k})}{|\vec{r} _{j}-%
\vec{r}_{k}|^{3}}=-\lambda (\vec{r}_{k}-\vec{r}{c})\qquad k=1,2,...,n,
\label{CCequation}
\end{equation}
where
\begin{equation}
\vec{c}=\frac{\vec{C}}{M_t}, \qquad \qquad \vec{C}=m_1\vec{r}_1+m_2\vec{r}%
_2+...+m_n\vec{r}_n,
\end{equation}
\begin{equation}
\lambda=\frac{U}{2I}, \qquad \qquad I=\frac{1}{2}\sum_{i=1}^n m_i \|\vec{r}%
_i\|^2.
\end{equation}
Due to the higher dimensions and degrees of freedom $n-$body problem has not been completely solved for $n>2$. Therefore a number of restriction techniques have been used to find special solutions of the few body problem. See for example \cite{shoaib5}, \cite{sekiguchi} and \cite{Rusu}. The two most common techniques used to reduce the dimension of the phase space are the consideration of symmetries are taking one of the masses to be infinitesimal. We consider four point masses on the vertices of an isosceles trapezoid with
two equal masses $m_1=m_4$ at positions $(\mp 0.5, r_B)$ and $m_2=m_3$ at
positions $(\mp \alpha/2, r_A)$. We derive, both analytically and
numerically, regions of central configurations in the phase space where it
is possible to choose positive masses.

The CC equations for a general 5-body problem derived from (\ref{CCequation}) are as under.
\begin{eqnarray}
\frac{{m_2}\vec{ r}_{12}}{\left| \vec{ r}_1-\vec{ r}_2\right| ^3}+\frac{{m_3}
\vec{ r}_{13}}{\left| \vec{ r}_1-\vec{ r}_3\right| ^3}+\frac{{m_4} \vec{ r}%
_{14}}{\left| \vec{ r}_1-\vec{ r}_4\right| ^3}=-\lambda \left(\vec{ r}_1-%
\vec{c}\right) \\
\frac{{m_1} \vec{ r}_{21}}{\left| \vec{ r}-\vec{ r}_2\right| ^3}+\frac{{m_3}
\vec{ r}_{23}}{\left| \vec{ r}_2-\vec{ r}_3\right| ^3}+\frac{{m_4} \vec{ r}%
_{24}}{\left| \vec{ r}_2-\vec{ r}_4\right| ^3}=-\lambda \left(\vec{ r}_2-%
\vec{c}\right) \\
\frac{{m_1}\vec{ r}_{31}}{\left| \vec{ r}_3-\vec{ r}_1\right| ^3}+\frac{{m_2}%
\vec{ r}_{32}}{\left| \vec{ r}_3-\vec{ r}_2\right| ^3}+\frac{{m_4} \vec{ r}%
_{34}}{\left| \vec{ r}_3-\vec{ r}_4\right| ^3}=-\lambda \left(\vec{ r}_3-%
\vec{ c}\right) \\
\frac{{m_1} \vec{ r}_{41}}{\left| \vec{ r}_4-\vec{ r}_1\right| ^3}+\frac{{m_2%
} \vec{ r}_{42}}{\left| \vec{ r}_4-\vec{ r}_2\right| ^3}+\frac{{m_3} \vec{ r}%
_{43}}{\left| \vec{ r}_4-\vec{ r}_3\right| ^3}=-\lambda \left(\vec{ r}_4-%
\vec{ c}\right)
\end{eqnarray}

\textbf{Theorem 1:}{\emph{\ Consider four bodies of masses $%
m_{1}=M,m_{2}=m=m_{3}$ and $m_{4}=M$. The four bodies are placed at the
vertices of a trapezoid
\begin{equation}
\mathbf{r}_{1}=(-0.5,-r_{B}),\mathbf{r}_{2}=(-\frac{\alpha }{2},r_{A}),%
\mathbf{r}_{3}=(\frac{\alpha }{2},r_{A})\text{ and }\mathbf{r}%
_{4}=(0.5,-r_{B}),  \label{coordinates}
\end{equation}%
shown in figure 1. Where $r_{A}$ is the distance from the centre of mass of
the system to the centre of mass of $m_{2}$ and $m_{3}$ and $r_{B}$ is the
distance from the centre of mass of the system to the centre of mass of $%
m_{1}$ and $m_{4}$. Then
\begin{eqnarray}
m = \frac{ab(a+b-2 a b ) }{(a + b) (a+b-2ab+\alpha(a-b))}  \label{m-general}
\\
M = \frac{a b \alpha( a-b )}{(a+b)(a+b-2ab+\alpha(a-b))}  \label{M-general}
\end{eqnarray}
where
\begin{equation}
a=\left( \left( 0.5-\frac{\alpha }{2}\right) {}^{2}+\beta {}^{2}\right)
^{3/2},{\ }b=\left( \left( 0.5+\frac{\alpha }{2}\right) {}^{2}+\beta
^{2}\right) ^{3/2}
\end{equation}
make the configuration $\mathbf{r}=(\vec{r}_1,\vec{r}_2,\vec{r}_3,\vec{r}_4)$
central.}}
\begin{figure}[tbp]
\centering
\resizebox{80mm}{!}{
\includegraphics {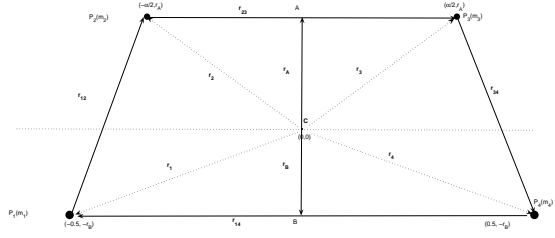}}
\caption{Trapezoidal Model}
\label{trapezoidalModel}
\end{figure}

\textbf{Theorem 2: }{\emph{ Let $\mathbf{r}=(\vec{r}_1,\vec{r}_2,\vec{r}_3,%
\vec{r}_4)$ be a central configuration as defined in theorem 1. Then there
exist a region $R={R_{f_1}}^c \cap {R_{f_3}}^c)$ in the $\alpha\beta$-plane such that for any $%
\alpha,\beta \in R$ there exists positive masses $(m,M,m,M)$ making $\mathbf{%
r}$ a central configuration. The regions $R_{f_1}$ and $R_{f_2}$ are
$$
R_{f_1}=\{(\alpha,\beta)|\alpha<g_1(\beta)\text{ and } 0<\beta<1\}
$$ $$
R_{f_3}=\{(\alpha,\beta)|\alpha<g_3(\beta)\text{ and } 0<\beta<1\}
$$ where
\begin{eqnarray}
g_1(\beta)=\frac{\sqrt{-2. \beta ^6-1.5 \beta ^4+2 \left(\beta ^2+0.25\right)^{3/2}-0.375 \beta ^2-0.03125}}{\sqrt{1.5 \beta ^4+\frac{-0.75 \beta ^2-0.375}{\sqrt{\beta ^2+0.25}}-0.09375}}.
\end{eqnarray}
$$
g_2(\beta)=\sqrt{-\frac{\sqrt{\text{h1}^2-4 \text{h0} \text{h2}}}{2 \text{h2}}-\frac{\text{h1}}{2 \text{h2}}}$$
Numerically, region $R$ is given by the colored part of figure (3b).
}}

Before we prove theorem 1, we recall a lemma given by Roy and Steves \cite{Steves}.

\textbf{Lemma 1:} {\emph{\ Let $\mathbf{r=r}_{A}-\mathbf{r}_{B}$ (ref:
figure 1) and then using the geometry of our proposed problem we arrive at
the following relationships between $\mathbf{r}_{i},i=1,2,3,4,$ $\mathbf{r}$
and $\mathbf{r}_{41}.$%
\begin{eqnarray*}
\mathbf{r}_{1} &=&-\frac{m}{M+m}\mathbf{r+}\frac{1}{2}\mathbf{r}_{41}, \\
\mathbf{r}_{2} &=&\frac{M}{m+M}\mathbf{r+}\frac{\alpha }{2}\mathbf{r}_{41},
\\
\mathbf{r}_{3} &=&\frac{M}{m+M}\mathbf{r-}\frac{\alpha }{2}\mathbf{r}_{41},
\\
\mathbf{r}_{4} &=&-\frac{m}{M+m}\mathbf{r-}\frac{1}{2}\mathbf{r}_{41}.
\end{eqnarray*}%
}}

\section{The proof of theorem 1}

Without loss of generality, it is assumed that the centre of mass of the
system $\vec{c}=0$, $\vec{r}_{23}=-\alpha \vec{r}_{41}$, $r_{BA}=|\vec{r}%
_{A}-\vec{r}_{B}|=\beta r_{41}$. This gives us
\begin{equation}
\vec{r}_B=-\frac{m}{M}\vec{r}_A \text{ and } r_A=\frac{M\beta}{m+M}
\end{equation}

Using these assumptions with equation (1) we obtain the following equations
of motion.
\begin{equation}
\ddot{\mathbf{r}}_{1}=\frac{m\mathbf{r}_{12}}{\left( \left( 0.5-\frac{\alpha
}{2}\right) ^{2}+\beta ^{2}\right) ^{3/2}}+\frac{m\mathbf{r}_{13}}{\left(
\left( 0.5+\frac{\alpha }{2}\right) ^{2}+\beta ^{2}\right) ^{3/2}}+M\mathbf{r%
}_{14},
\end{equation}%
\begin{equation}
\ddot{\mathbf{r}}_{2}=\frac{M\mathbf{r}_{21}}{\left( \left( 0.5-\frac{\alpha
}{2}\right) ^{2}+\beta ^{2}\right) ^{3/2}}+\frac{m\mathbf{r}_{23}}{\alpha
^{3}}+\frac{M\mathbf{r}_{24}}{\left( \left( 0.5+\frac{\alpha }{2}\right)
^{2}+\beta ^{2}\right) ^{3/2}},
\end{equation}%
\begin{equation}
\ddot{\mathbf{r}}_{3}=\frac{M\mathbf{r}_{31}}{\left( \left( 0.5+\frac{\alpha
}{2}\right) ^{2}+\beta ^{2}\right) ^{3/2}}+\frac{m\mathbf{r}_{32}}{\alpha
^{3}}+\frac{M\mathbf{r}_{34}}{\left( \left( 0.5-\frac{\alpha }{2}\right)
^{2}+\beta ^{2}\right) ^{3/2}},
\end{equation}%
\begin{equation}
\ddot{\mathbf{r}}_{4}=\frac{m\mathbf{r}_{42}}{\left( \left( 0.5+\frac{\alpha
}{2}\right) ^{2}+\beta ^{2}\right) ^{3/2}}\ +M\mathbf{r}_{41}+\frac{m\mathbf{%
r}_{43}}{\left( \left( 0.5-\frac{\alpha }{2}\right) ^{2}+\beta ^{2}\right)
^{3/2}}.
\end{equation}

It is clear from lemma 1 that it is enough to study the equations for $%
\mathbf{r=}\frac{m+M}{2M}\left( \mathbf{r}_{2}+\mathbf{r}_{3}\right) $ and $%
\mathbf{r}_{41}$ as $\mathbf{r}_{i},i=1,2,3,4,$ are linear combination of $%
\mathbf{r}$ and $\mathbf{r}_{41}.$%
\begin{equation}
\ddot{\mathbf{r}}_{41}=[-2M+\frac{m}{a}\left( \alpha -1\right) -\frac{m}{b}%
(\alpha +1)]\mathbf{r}_{41}.{\ \ \ \ \ \ \ \ \ \ \ \ }  \label{r41}
\end{equation}%
\begin{equation}
\ddot{\mathbf{r}}_{3p2}=-[\frac{m+M}{a}+\frac{m+M}{b}]\left( \mathbf{r}_{2}+%
\mathbf{r}_{3}\right) .{\ \ \ \ \ \ \ \ \ \ \ \ \ }  \label{r32}
\end{equation}%
Using equation (\ref{CCequation}) in conjunction with equations (\ref{r41})
and (\ref{r32}) we obtain the following equations of central configurations
for the trapezoidal four-body problem.

\begin{equation}
2M-\frac{m(\alpha -1)}{a}+\frac{m(\alpha +1)}{b}=\lambda ,  \label{CCeq1}
\end{equation}%
\begin{equation}
\frac{m+M}{a}+\frac{m+M}{b}=\lambda .  \label{CCeq2}
\end{equation}
It is a straightforward exercise to solve (\ref{CCeq1}) and (\ref{CCeq2}) to
obtain (\ref{m-general}) and (\ref{M-general}). This completes the proof of
theorem 1.
\begin{figure}[tbp]
\centering
\resizebox{130mm}{!}{
\includegraphics {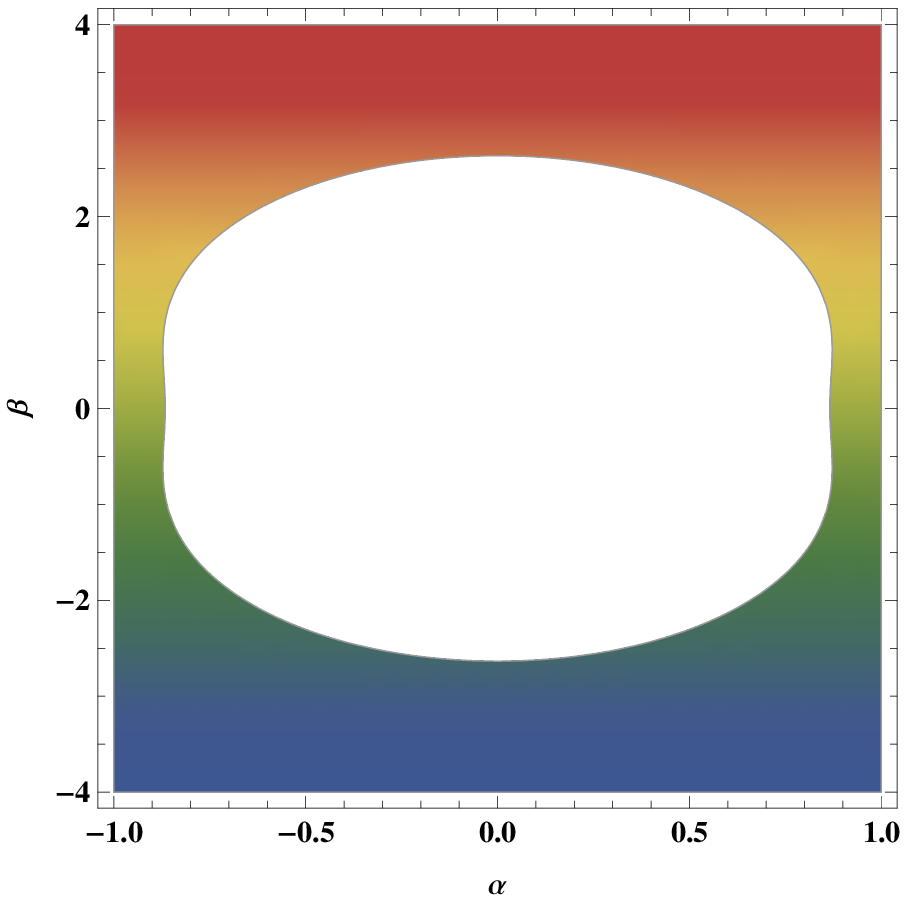}
\includegraphics {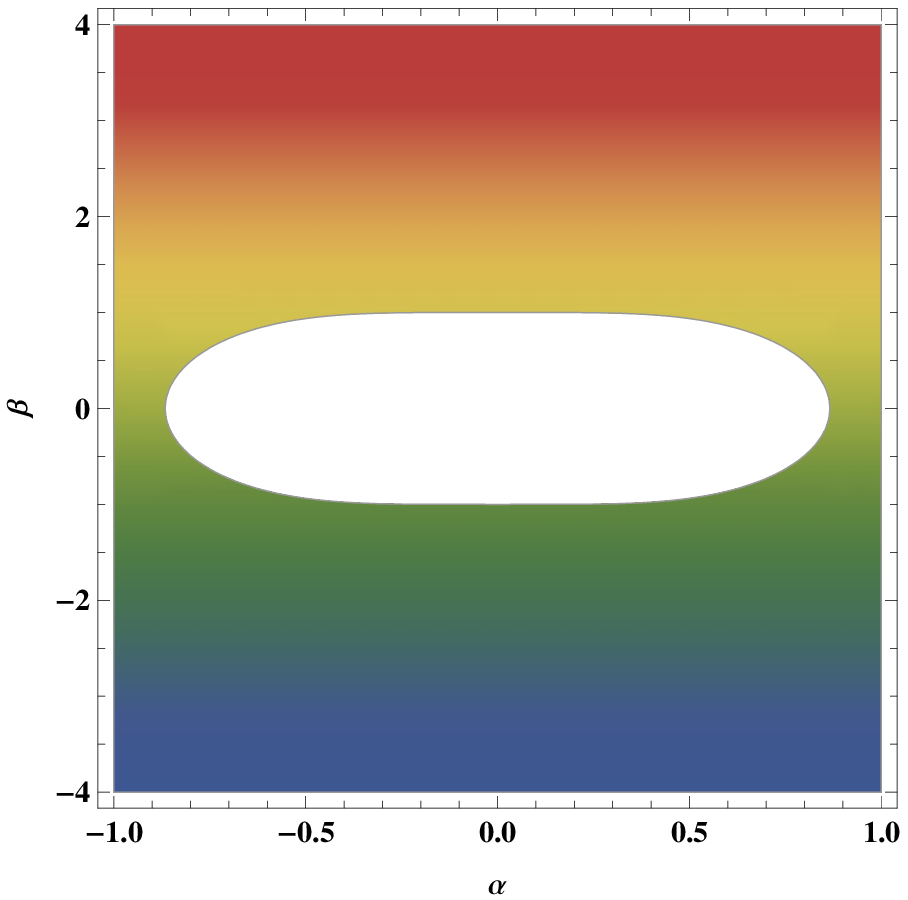}}
\caption{\hskip0.5cm a ) $R_{f_1}$ (white) where $f_1>0$. \hskip1cmb)  $R_{f_3}$ (white) where $f_3>0$ }
\label{Rf1and3}
\end{figure}
\begin{figure}[tbp]
\centering
\resizebox{130mm}{!}{
\includegraphics {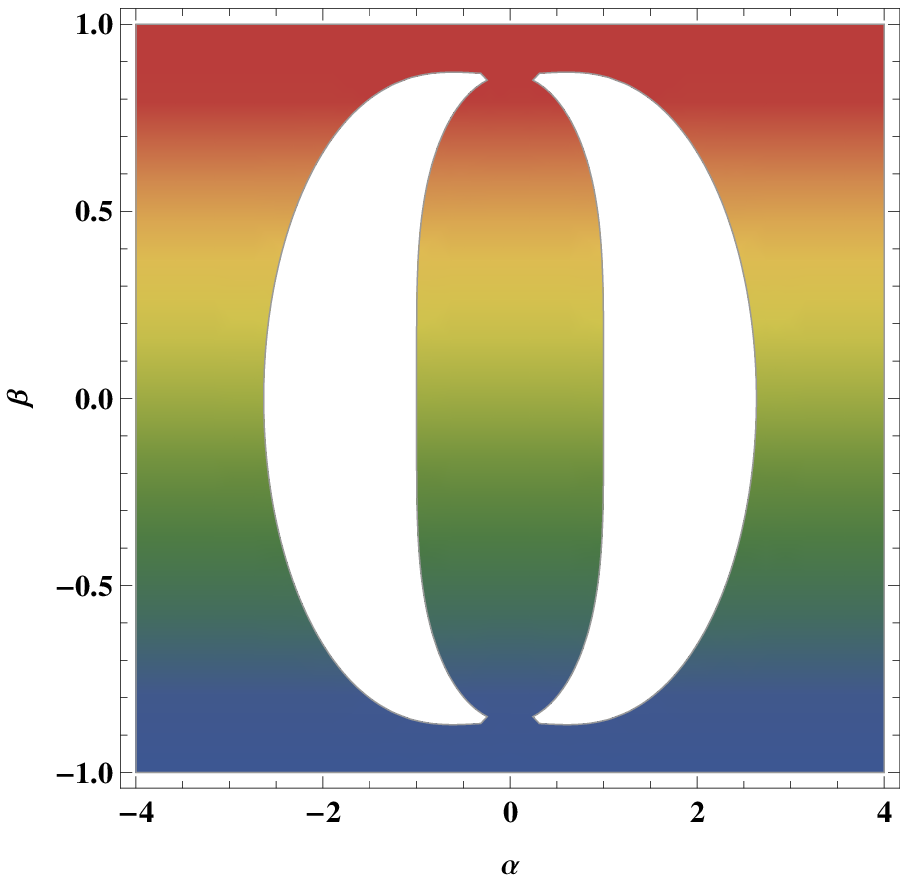}
\includegraphics {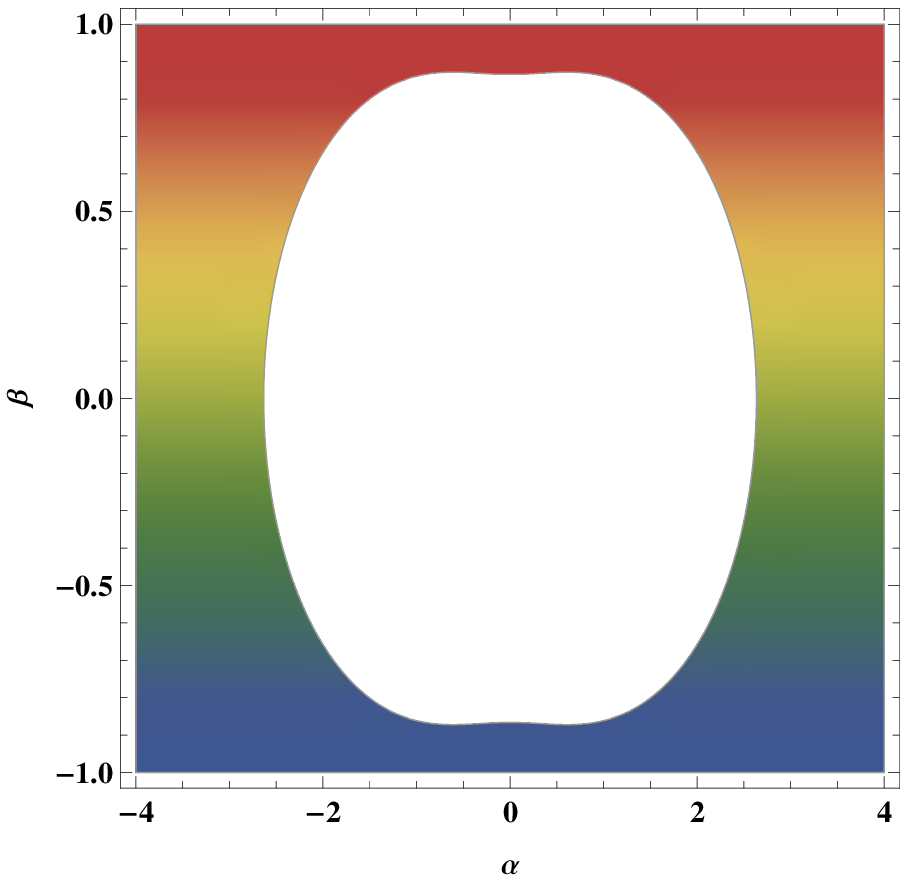}}
\caption{\hskip0.5cm  a) Region (colored), $R_m=(R_{f_1}\cap R_{f_3})\cup ({R_{f_1}}^c \cap {R_{f_3}}^c)$  b) The central configuration region ($R$) where both $m$ and $M$ are positive ( $R={R_{f_1}}^c \cap {R_{f_3}}^c)$).}
\label{Rm}
\end{figure}

\section{Proof of theorem 2}

Let
\begin{equation*}
f_{1}=a+b-2ab,f_{2}=a-b,f3=a+b-2ab+\alpha (a-b)
\end{equation*}%
To find the region where the mass function $m$ is positive we need to find
regions where

\begin{enumerate}
\item $f_1>0$,$f_3>0$,

\item $f_1<0$,$f_3<0$.
\end{enumerate}

Similarly, for the mass function $M$ to be positive

\begin{enumerate}
\item $f_2>0$,$f_3>0$

\item $f_2<0$,$f_3<0$
\end{enumerate}

To do a sign analysis of $f_{1}$, we will need to solve $f_{1}=0$. Its nearly impossible to solve $f_1>0$ . Therefore we write its polynomial approximation.
\begin{eqnarray}
f_1aprox&=&\alpha ^2 \left(-1.5 \beta ^4+\frac{0.75 \beta ^2}{\sqrt{\beta ^2+0.25}}+\frac{0.375}{\sqrt{\beta ^2+0.25}}+0.09375\right)-2 \beta ^6-1.5 \beta ^4\nonumber\\&+&2 \sqrt{\beta ^2+0.25} \beta ^2-0.375 \beta ^2+0.5 \sqrt{\beta ^2+0.25}-0.03125.\end{eqnarray}
Now it is a straightforward exercise to show that $f_1>0$ in $R_{f_1}$.
$$
R_{f_1}=\{(\alpha,\beta)|\alpha<g_1(\beta)\text{ and } 0<\beta<1\}
$$ where
\begin{eqnarray}
g_1(\beta)=\frac{\sqrt{-2. \beta ^6-1.5 \beta ^4+2 \left(\beta ^2+0.25\right)^{3/2}-0.375 \beta ^2-0.03125}}{\sqrt{1.5 \beta ^4+\frac{-0.75 \beta ^2-0.375}{\sqrt{\beta ^2+0.25}}-0.09375}}.
\end{eqnarray}
Numerically, $R_{f_1}$ is given in figure (\ref{Rf1and3}a).  The common denominator of $m$, and $M$  can be analyzed in a similar way.
\begin{eqnarray}
f_3aprox&=&h_2\alpha ^4+h_1\alpha ^2+h_0\end{eqnarray}
where
\begin{eqnarray}
h_0&=&-2 \beta ^6-1.5 \beta ^4+2 \left(\beta ^2+0.25\right)^{3/2}-0.38 \beta ^2-0.03\nonumber\\
h_1&=&-1.5 \beta ^4-\frac{0.75 \beta ^2}{\sqrt{\beta ^2+0.25}}+0.1\nonumber\\
h_2&=&(-0.375 \beta ^6-0.09 \beta ^4+\left(-0.19 \sqrt{\beta ^2+0.25}-0.023\right) \beta ^2-0.031 \sqrt{\beta ^2+0.25}\nonumber\\&+&\frac{0.047 \beta ^4}{\sqrt{\beta ^2+0.25}}-0.006)\frac{1}{\left(1. \beta ^2+0.25\right)^2}
\end{eqnarray}
Using  approximate techniques with help from symbolic computation in Mathematica it can be shown that $f_3$ attains positive values in the following region.

$$
R_{f_3}=\{(\alpha,\beta)|\alpha<g_3(\beta)\text{ and } 0<\beta<1\}
$$ where
$$
g_2(\beta)=\sqrt{-\frac{\sqrt{\text{h1}^2-4 \text{h0} \text{h2}}}{2 \text{h2}}-\frac{\text{h1}}{2 \text{h2}}}$$
Numerically, $R_{f_3}$ is given in figure (\ref{Rf1and3}b).

Therefore the central configuration region where $m>0$ is given by
$$
R_m=(R_{f_1}\cap R_{f_3})\cup ({R_{f_1}}^c \cap {R_{f_3}}^c)
$$
Numerically, $R_m$ is given in figure (\ref{Rm}).  As $f_2<0$ for all values of $\alpha$ and $\beta$ therefore $M>0$ in the region where $f_3<0$. This region is numerically represented by the colored part of figure (\ref{Rf1and3}b) and analytically by the compliment of ${R_{f_3}}.$
Hence the central configuration region ($R$) where both $m$ and $M$ are positive is given by $R=({R_{f_1}}^c \cap {R_{f_3}}^c)$. Numerically this region is given in figure (\ref{Rm}b).
\section{Conclusions}
In this paper we model  non-collinear trapezoidal four-body problem where the masses are
placed at the vertices of an isosceles trapezoid. Expressions for $m$ and $M$  are formed as functions of $\alpha$, and  $\beta$ which gives central configurations in the trapezoidal four-body problems. We show that in the $\alpha \beta$-plane, $m$ is positive when $(\alpha,  \beta) \in R_m$. Similarly $M$ is positive when  $(\alpha,  \beta) \in R_{f_3}^c$. We have identified regions in the $\alpha\beta$-plane where no central configurations are possible. A central configuration region  $R=({R_{f_1}}^c \cap {R_{f_3}}^c)$ for the isosceles trapezoidal 4-body problem is identified in the $\alpha \beta-$ plane  where $m$ and $M$ are both  positive. No central configurations are possible outside this region unless we allow one of the masses to become negative.

\vskip1cm
{\bf{Acknowledgement:}} The author thanks the Deanship of research at the
University of Hail, Saudi Arabia for funding this work under grant
number SM14014.

\end{document}